\documentclass[a4paper,12pt]{article}
\usepackage{amsfonts,amssymb,graphicx}

\newtheorem{theorem}{Theorem}

\begin{document}

\author{
Octavian G. Mustafa\footnote{On leave from the University of Craiova, Faculty of Mathematics and Computer Science, Al{.} I{.} Cuza 13, Craiova, Romania.}\\
\small{University of Vienna, Faculty of Mathematics,}\\
\small{Nordbergstrasse 15, A-1090 Vienna, Austria}\\
\small{e-mail address: octaviangenghiz@yahoo.com}\\
\small{Fax number: +40251525758}\\
and\\
Donal O'Regan\\
\small{National University of Ireland,}\\
\small{Mathematics Department,}\\
\small{Galway, Ireland}\\
\small{e-mail address: donal.oregan@nuigalway.ie}
}

\title{On the Nagumo uniqueness theorem}
\date{}
\maketitle

\noindent{\bf Abstract} By a convenient reparametrisation of the integral curves of a nonlinear ordinary differential equation (ODE), we are able to improve the conclusions of the recent contribution [A. Constantin, Proc. Japan Acad. {\bf 86(A)} (2010), 41--44]. In this way, we establish a flexible uniqueness criterion for ODEs without Lipschitz-like nonlinearities.

\noindent{\bf Key-words:} Ordinary differential equation; uniqueness of solution; repara\-metrization

\noindent{\bf Classification:} 34A12; 34A34

\section{Introduction}

Let us consider the nonlinear ordinary differential equation
\begin{eqnarray}
x^{\prime}+f(t,x)=0,\quad t>0, \label{main_eq}
\end{eqnarray}
where the nonlinearity $f$ is assumed continuous everywhere and $f(t,0)\equiv0$. In this note, we shall discuss a set of very permissive conditions to be imposed on $f$ so that the only solution of (\ref{main_eq}) starting from $x(0)=0$ will be the trivial solution.

Several results in the recent literature use such uniqueness criteria to perform a phase plane analysis of various differential equations and systems from the applied sciences, see e{.}g{.}, \cite{constdyn,balabane,kk}.

One of the most powerful uniqueness theorems regarding (\ref{main_eq}) is due to Nagumo \cite{nag} and states that, if $\lim\limits_{t\searrow0}f(t,x)=0$ uniformly with respect to $x$ in $[-1,1]$ and
\begin{eqnarray}
\vert f(t,x_{1})-f(t,x_{2})\vert\leq\frac{1}{t}\cdot\vert x_{1}-x_{2}\vert,\quad t>0,\thinspace \vert x_{1}\vert,\thinspace\vert x_{2}\vert\leq1,\label{nagu}
\end{eqnarray}
the initial value problem associated to the equation possesses a unique solution. The result has been generalized to $n$--th order equations in \cite{nag2,wintner,c97} and other papers.

In a different direction \cite{ath}, given the smooth function $u:[0,1]\rightarrow[0,+\infty)$ with $u(0)=0$ and $u^{\prime}(t)>0$ everywhere in $(0,1]$, if the nonlinearity $f$ verifies the restriction
\begin{eqnarray}
\vert f(t,x_{1})-f(t,x_{2})\vert\leq\frac{u^{\prime}(t)}{u(t)}\cdot\vert x_{1}-x_{2}\vert,\quad t\in(0,1],\thinspace \vert x_{1}\vert,\thinspace\vert x_{2}\vert\leq1,\label{atha}
\end{eqnarray}
and $\lim\limits_{t\searrow0}\frac{f(t,x)}{u^{\prime}(t)}=0$ uniformly with respect to $x$ in $[-1,1]$ then the only solution of (\ref{main_eq}) with null initial datum will be the null solution.

It has been established in \cite{cjapan} that the restriction (\ref{atha}) can be replaced with
\begin{eqnarray}
\vert f(t,x)\vert\leq\frac{u^{\prime}(t)}{u(t)}\cdot\omega(\vert x\vert),\quad t\in(0,1],\thinspace \vert x\vert\leq1,\label{adri}
\end{eqnarray}
where the comparison function $\omega:[0,1]\rightarrow[0,+\infty)$ is continuous and increasing, $\omega(0)=0$ and $\omega(r)>0$ for $r>0$, and also
\begin{eqnarray}
\int_{0}^{r}\frac{\omega(s)}{s}ds\leq r,\quad r\in[0,1].\label{hyp_omega}
\end{eqnarray}
For a variant of this result in the case of second order ODEs and other generalizations, see \cite{Mustafa,Mejstrik}.

In the following sections we present an extension of the result from \cite{cjapan} and comment on the connections between the theorems of Nagumo, Athanassov and Constantin.

\section{The result}

Assume that there exist two smooth functions $v$, $\lambda:[0,1]\rightarrow[0,+\infty)$ such that $v(0)=\lambda(0)=0$ and $v^{\prime}(t)$, $\lambda(t)>0$ throughout $(0,1]$. Consider also that
\begin{eqnarray*}
\frac{\lambda(t)f(t,x)}{v(t)},\thinspace\frac{f(t,x)}{v^{\prime}(t)}\longrightarrow0\mbox{ as }t\searrow0\mbox{ uniformly in }\vert x\vert\leq 1.
\end{eqnarray*}
As a consequence, we presume that
\begin{eqnarray*}
\vert\lambda(t)f(t,x)\vert\leq v(t),\quad t\in[0,1],\thinspace\vert x\vert\leq1.
\end{eqnarray*}

We have the following result.

\begin{theorem}\label{theo_m}
Let $\int_{0+}^{1}\frac{v(w)}{\lambda(w)}dw<+\infty$. Given a continuous, strictly increasing function $\omega:[0,1]\rightarrow[0,+\infty)$, with $\omega(0)=0$, such that
\begin{eqnarray}
\int_{0+}^{t}\omega(\varepsilon v(w))\frac{dw}{\lambda(w)}\leq\varepsilon v(t)\mbox{ for all }\varepsilon>0\mbox{ and }t\in(0,1],\label{main_ineq_2}
\end{eqnarray}
suppose further that, for $t\in(0,1]$ and $\vert x\vert\leq 1$, one has
\begin{eqnarray}
\vert f(t,x)\vert\leq\frac{\omega(\vert x\vert)}{\lambda(t)}.\label{main_ineq_1}
\end{eqnarray}
Then the only solution of (\ref{main_eq}) starting from $x(0)=0$ is the trivial solution.
\end{theorem}

{\bf Proof.} For some $T\in(0,1)$, introduce the function $\tau:(0,T]\rightarrow[\tau_{-},\tau_{+})$, with $-\infty<\tau_{-}<\tau_{+}\leq+\infty$, by means of
\begin{eqnarray}
\left\{
\begin{array}{ll}
\tau(t)=\tau_{-}+\int_{t}^{T}\frac{ds}{\lambda(s)},\quad 0<t\leq T,\\
\\
\tau(T)=\tau_{-},\quad\lim\limits_{t\searrow0}\tau(t)=\tau_{+}.
\end{array}
\right.\label{t_tau_0}
\end{eqnarray}

We assume for the sake of contradiction that the initial value problem associated to the equation (\ref{main_eq}) has a non-trivial solution $x$ with $\vert x(t)\vert\leq1$ for $t\in[0,T]$.

The proof relies on a change of variables, that is, we introduce the function $y$ via  $y(\tau(t))=x(t)$. Now,
\begin{eqnarray*}
\frac{dy}{d\tau}&=&\frac{\frac{dx}{dt}}{\frac{d\tau}{dt}}=-\lambda(t(\tau))\cdot f(t(\tau),x(t(\tau)))=-\lambda(t(\tau))f(t(\tau),y(\tau))\\
&=&F(\tau,y),
\end{eqnarray*}
where $\vert F(\tau,y)\vert\leq v(t)$ and the quantity $t(\tau)$ is given by
\begin{eqnarray}
\left\{
\begin{array}{ll}
t(\tau)=\int_{\tau}^{\tau_{+}}\lambda(t(s))ds,\quad \tau_{-}\leq \tau<\tau_{+},\\
\\
t(\tau_{-})=T=\int_{\tau_{-}}^{\tau_{+}}\lambda(t(s))ds<+\infty,\quad\lim\limits_{\tau\nearrow\tau_{+}}t(\tau)=0.
\end{array}
\right.\label{t_tau_1}
\end{eqnarray}

Further, we introduce the continuous function $\alpha$ from $v(t)=\alpha(\tau(t))$, $t\in(0,T]$, or, equivalently, $\alpha(\tau)=v(t(\tau))$. We claim that the function is a member of $L^{1}((\tau_{-},\tau_{+}),\mathbb{R})$. This follows from
\begin{eqnarray*}
\int_{\tau}^{\tau_{+}}\alpha(\xi)d\xi&=&\int_{t}^{0}\alpha(\tau(w))\cdot\frac{d\tau}{dw}\thinspace dw=\int_{0}^{t}\frac{\alpha(\tau(w))}{\lambda(w)}dw\\
&=&\int_{0+}^{t}\frac{v(w)}{\lambda(w)}dw<+\infty.
\end{eqnarray*}

Since $\lim\limits_{\tau\nearrow\tau_{+}}y(\tau)=\lim\limits_{t\searrow0}x(t)=0$, the integrability of $\alpha$ leads to
\begin{eqnarray*}
-y(\tau)=\lim\limits_{\tau\rightarrow\tau_{+}}y(\tau)-y(\tau)=-\int_{\tau}^{\tau_{+}}F(\xi,y(\xi))d\xi
\end{eqnarray*}
and respectively
\begin{eqnarray*}
\vert y(\tau)\vert&\leq&\int_{\tau}^{\tau_{+}}\vert F(\xi,y(\xi))\vert d\xi\leq\int_{\tau}^{\tau_{+}}\omega(\vert y(\xi)\vert)d\xi\\
&<&\int_{\tau}^{\tau_{+}}\omega\left(\sup\limits_{s\in[\tau,\tau_{+})}\frac{\vert y(s)\vert}{\alpha(s)}\cdot\alpha(\xi)\right)d\xi\\
&\leq&\sup\limits_{s\in[\tau,\tau_{+})}\frac{\vert y(s)\vert}{\alpha(s)}\cdot\alpha(\tau).
\end{eqnarray*}
The latter estimate follows from (\ref{main_ineq_2}) for $\varepsilon=\sup\limits_{s\in[\tau(t),\tau_{+})}\frac{\vert y(s)\vert}{\alpha(s)}$ and $t\in(0,T]$ fixed.

In conclusion,
\begin{eqnarray}
\frac{\vert y(\tau)\vert}{\alpha(\tau)}<\sup\limits_{s\in[\tau,\tau_{+})}\frac{\vert y(s)\vert}{\alpha(s)},\quad\tau\in(\tau_{-},\tau_{+}).\label{fin_est_nagg}
\end{eqnarray}
Since $\lim\limits_{\tau\nearrow\tau_{+}}\frac{y(\tau)}{\alpha(\tau)}=\lim\limits_{t\searrow0}\frac{x(t)}{v(t)}=0$, taking into account the continuity of the functions involved, we deduce that $\sup\limits_{s\in[\tau,\tau_{+})}\frac{\vert y(s)\vert}{\alpha(s)}=\max\limits_{s\in[\tau,\tau_{+})}\frac{\vert y(s)\vert}{\alpha(s)}=\frac{\vert y(\xi_{\tau})\vert}{\alpha(\xi_{\tau})}$ for some $\xi_{\tau}\in[\tau,\tau_{+})$. This leads to a contradiction if we take $\tau=\xi_{\tau}$ in (\ref{fin_est_nagg}).

The proof is complete. $\square$

\section{Discussion}

First of all, notice that for (\ref{adri}), one has $\lambda=\frac{u}{u^{\prime}}$ and $v=u$. Thus, (\ref{main_ineq_2}) reduces to (\ref{hyp_omega}), see also \cite[p. 42]{cjapan}. As a consequence, Theorem \ref{theo_m} is a generalization of the recent result of Constantin.

We would like to offer here a clarification of the technical change of variables from \cite{cjapan}. To this end, let us remark that, at least theoretically, {\it the Athanassov theorem from \cite{ath} is equivalent to the classical Nagumo uniqueness criterion\/}. Our claim is supported by the next computations. Take $y(u(t))=x(t)$ and observe that
\begin{eqnarray*}
\frac{dy}{du}=\frac{\frac{dx}{dt}}{\frac{du}{dt}}=\frac{f(t,y(u))}{u^{\prime}(t)}=\frac{dt}{du}\cdot f(t(u),y(u))=g(u,y(u))
\end{eqnarray*}
and
\begin{eqnarray*}
\vert g(u,y(u))\vert=\frac{\vert f(t,y)\vert}{u^{\prime}(t)}\leq\frac{\vert y\vert}{u},\quad u>0.
\end{eqnarray*}
These restrictions match precisely the Nagumo hypotheses in \cite{nag}. It is obvious that, in practice, inverting the general function $u=u(t)$ is not an easy task. Consequently, the previous equivalence remains just a remark.

So, one might wonder, if the Nagumo and Athanassov conditions are almost the same, maybe there are other changes of variables that will lead in various circumstances to some sharp uniqueness criterion. One of these changes of variables has been proposed very recently in \cite{cjapan}. In fact, the change of variables reads as follows
\begin{eqnarray}
u(t)=u(t(\tau))=c\cdot\mbox{e}^{-\tau},\quad\tau_{-}\leq\tau<+\infty=\tau_{+},\label{ch_adr}
\end{eqnarray}
for some constant $c>0$.

We have that
\begin{eqnarray*}
\frac{d}{d\tau}[u(t(\tau))]&=&-u(t(\tau))=-c\cdot\mbox{e}^{-\tau}\\
&=&\frac{du}{dt}\cdot\frac{dt}{d\tau},
\end{eqnarray*}
leading to
\begin{eqnarray*}
\frac{dt}{d\tau}=-\frac{u}{u^{\prime}}=-\lambda(t(\tau)),\quad t(\tau)=\int_{\tau}^{+\infty}\lambda(t(s))ds.
\end{eqnarray*}
The latter integral equation is the same as (\ref{t_tau_1}). The inverse of its solution $t=t(\tau)$, which is nothing but a recast of (\ref{t_tau_0}), follows naturally from a {\it reparametrization\/} of the graph of function $u$ from (\ref{atha}). We have (recall that $u(0)=0$)
\begin{eqnarray*}
u(t)&=&\int_{0}^{t}u^{\prime}(s)ds=\int_{\tau_{+}}^{\tau(t)}u^{\prime}(s(\xi))\cdot\frac{ds}{d\xi}\thinspace d\xi=\int_{\tau(t)}^{\tau_{+}}u^{\prime}(s(\xi))\lambda(s(\xi))d\xi\\
&=&\int_{\tau(t)}^{+\infty}u(s(\xi))d\xi.
\end{eqnarray*}

Such reparametrizations have been employed successfully in various problems regarding ordinary and partial differential equations, see e{.}g{.}, \cite[p. 17]{chicone} or \cite[p. 134]{musta}.

Finally, let us consider a more general reparametrization $t\mapsto\tau\mapsto t$ given by the differential equations
\begin{eqnarray*}
\frac{dt}{d\tau}=-\lambda(t(\tau),\tau),\quad\frac{d\tau}{dt}=-\frac{1}{\lambda(t,\tau(t))},
\end{eqnarray*}
and, following (\ref{ch_adr}), put $u(\tau,t(\tau))=c\cdot\mbox{e}^{-\tau}$.

Further, we have
\begin{eqnarray*}
\frac{d}{d\tau}[u(\tau,t(\tau))]&=&\frac{\partial u}{\partial\tau}+\frac{\partial u}{\partial t}\cdot\frac{d t}{d\tau}=\frac{\partial u}{\partial\tau}+\frac{\partial u}{\partial t}\cdot[-\lambda(t(\tau),\tau)]\\
&=&-u(\tau,t(\tau)),
\end{eqnarray*}
yielding
\begin{eqnarray*}
\lambda(t,\tau)=\frac{u+\frac{\partial u}{\partial \tau}}{\frac{\partial u}{\partial t}},\quad u=u(\tau,t).
\end{eqnarray*}

Taking $u(\tau,t)=u(t)+\frac{1}{\tau}=c\cdot\mbox{e}^{-\tau}$, where $u$ satisfies the properties introduced by Athanassov \cite{ath}, we have either $\tau_{+}=+\infty$ or $\tau_{+}\mbox{e}^{\tau_{+}}=\frac{1}{c}$ and extract $t$ such that
\begin{eqnarray}
u(t(\tau))=c\cdot\mbox{e}^{-\tau}-\frac{1}{\tau}.\label{fin_gen_eq}
\end{eqnarray}

The key restriction (\ref{adri}) would be relaxed in this case up to
\begin{eqnarray*}
\vert f(t,x)\vert\leq\frac{u^{\prime}(t)}{u(t)-\frac{1}{[\tau(t)]^{2}}}\cdot\omega(\vert x\vert),\quad t\in(0,T],
\end{eqnarray*}
where $\tau=\tau(t)$ is the inverse of $t=t(\tau)$ from (\ref{fin_gen_eq}). It is still unclear to us how/if this will lead to a drastic improvement of the conclusions from Theorem \ref{theo_m} and \cite{cjapan}.

\small{

}
\end{document}